\documentclass{article}
\usepackage[english]{babel}
\usepackage{latexsym,amssymb}
\usepackage[cp1251]{inputenc}
\usepackage{amsfonts,amssymb}
\usepackage{euscript}
\newcounter{theorem} %[section]
\newcounter{lemma} %[section]

\begin{document}

\vspace*{7mm}

\Large

%\addtocounter{page}{-1}
%\renewcommand{\baselinestretch}{1.35}
\thispagestyle{empty}

%\noindent УДК 517.548\\
\Large
 \begin{center}
    \textbf{Vitalii Shpakivskyi\footnote{Institute of Mathematics of the NAS of Ukraine, Kyiv, Ukraine; shpakivskyi86@gmail.com}}
 \end{center}
%\vspace{3mm}
\vspace{7mm}
\Huge
 \begin{center}
   \textbf{
PART 2}
 \end{center}
\vspace{3mm}
%\vspace{7mm}

\LARGE
 \begin{center}
   \textbf{Algebras with Variable Structural Constants and Their Applications}
 \end{center}
\vspace{5mm}

% У роботі розвинуто метод П.~Гумберта  розв'язування диференціальних рівнянь з частинними похідними. Даний метод застосовано до
% побудови розв'язків тривимірних рівнянь Лапласа, Гельмгольца і Пуассона в класах голоморфних функцій.

%\vspace{15mm}

%\noindent {\bf  Vitalii Shpakivskyi}
%\vspace{3mm}

%\noindent {\bf Holomorphic functions of several hypercomplex variables and their applications}

%\vspace{7mm}

% The work developed P.~Humbert's method of solving PDEs. This method is applied to
% construction of solutions of the three-dimensional equations of Laplace, of Helmholtz and of Poisson in the classes of holomorphic functions.
 
\Large
%\newpage 
%\section{Основні поняття}

\section{Algebras with Variable Structural Constants}

Let $\mathbb{K}\in\{\mathbb{R},\mathbb{C}\}$ and let $D$ be a domain in $\mathbb{K}^m$.
Let $A$ be an $n$-dimensional associative algebra over the field $\mathbb{K}$ with basis $\{e_1,e_2,\ldots,e_n\}$.
%In this case, $1\leq m \leq n$ if $\mathbb{K}=\mathbb{R}$ and $1\leq m \leq 2n$ if $\mathbb{K}=\mathbb{C}$.
The multiplication table of the algebra $A$ has the form
$$
e_ie_j=\sum_{k=1}^{n}\gamma_{ij}^{k}\,e_k,
\qquad i,j=1,\ldots,n,
$$
where $\gamma_{ij}^{k}\in\mathbb{K}$.

To each point $x\in D$, we associate an algebra $A_x$ according to the rule
\begin{equation}\label{tabl-mnozh-algebra-zi-zminnymy-funkciamy}
e_ie_j=\sum_{k=1}^{n}\gamma_{ij}^{k}(x)\,e_k,
\qquad i,j=1,\ldots,n,
\end{equation}
where $\gamma_{ij}^{k}:D\rightarrow\mathbb{K}$ are given functions.

Thus, for a fixed collection of structural functions
$$
\gamma:=\{\gamma_{ij}^{k}:D\rightarrow\mathbb{K},\,\, i,j,k=1,\ldots,n\},
$$
each point $x$ of the domain $D$ is associated with an algebra $A_x(\gamma)$.
The family of algebras
$$
\mathcal{A}(D,\gamma):=\{A_x(\gamma)\}_{x\in D}
$$
will be called an \textbf{algebra with variable structural constants} $\gamma$.

\vskip2mm

This definition generalizes the classical definition of an algebra \cite{Pierce}.

\textbf{Remark 1.}
If all structural functions $\gamma_{ij}^{k}(x)$ are constant on the domain $D$, then all algebras $A_x(\gamma)$ coincide, and the algebra $\mathcal{A}(D,\gamma)$ reduces to an ordinary finite-dimensional algebra over the field $\mathbb{K}$.

Let us consider some examples of algebras with variable structural constants $\mathcal{A}(D,\gamma)$ (families of algebras).

\vskip2mm
\textbf{Example 1.}
Let $a,b:D\to\mathbb{K}$ be given functions. Consider a family of two-dimensional algebras with basis $\{e_1,e_2\}$ whose multiplication table has the form
\begin{equation}\label{tabl-mnozh-2-vym-algebra-zi-zminnymy-funkciamy}
\begin{array}{c||c|c|}
\cdot & e_1 & e_2\\
\hline\hline
e_1 & e_1 & e_2\\
\hline
e_2 & e_2 & a(x)e_1+b(x)e_2\\
\hline
\end{array}\,.
\end{equation}

For each fixed $x\in D$, we obtain a commutative associative algebra $A_x$.
If the functions $a$ and $b$ are constant, then this family reduces to the classical two-dimensional commutative algebra with constant structural constants.

\vskip2mm
\textbf{Example 2.}
Let $\lambda:D\to\mathbb{K}$ be a given function. Consider a family of two-dimensional algebras with basis $\{e_1,e_2\}$ whose multiplication table has the form
$$
\begin{array}{c||c|c|}
\cdot & e_1 & e_2\\
\hline\hline
e_1 & e_1 & e_2\\
\hline
e_2 & \lambda(x)e_2 & 0\\
\hline
\end{array}\,.
$$

If $\lambda(x)\neq 1$, then
$$
e_1e_2=e_2\neq \lambda(x)e_2=e_2e_1,
$$
and hence the algebra $A_x$ is noncommutative.

\vskip2mm
\textbf{Example 3.}
Let $a:D\to\mathbb{K}$ be a given function. Consider a family of three-dimensional algebras with basis $\{e_1,e_2,e_3\}$ whose multiplication table has the form
$$
\begin{array}{c||c|c|c|}
\cdot & e_1 & e_2 & e_3\\
\hline\hline
e_1 & 0 & a(x)e_3 & 0\\
\hline
e_2 & -a(x)e_3 & 0 & 0\\
\hline
e_3 & 0 & 0 & 0\\
\hline
\end{array}\,.
$$

For each $x\in D$, we obtain a noncommutative algebra. Indeed,
$$
e_1e_2=a(x)e_3,
\qquad
e_2e_1=-a(x)e_3,
$$
and if $a(x)\neq 0$, then $e_1e_2\neq e_2e_1$.

\vskip2mm
\textbf{Example 4.}
Let $\alpha,\beta:D\to\mathbb{K}$ be given functions. Consider a family of four-dimensional algebras with basis
$\{1,e_1,e_2,e_3\}$ whose multiplication table has the form
$$
\begin{array}{c||c|c|c|c|}
\cdot & 1 & e_1 & e_2 & e_3\\
\hline\hline
1   & 1 & e_1 & e_2 & e_3\\
\hline
e_1 & e_1 & -\alpha(x) & e_3 & -\alpha(x)e_2\\
\hline
e_2 & e_2 & -e_3 & -\beta(x) & \beta(x)e_1\\
\hline
e_3 & e_3 & \alpha(x)e_2 & -\beta(x)e_1 & -\alpha(x)\beta(x)\\
\hline
\end{array}\,.
$$

For each fixed $x\in D$, we obtain a generalized quaternion algebra.
If the functions $\alpha$ and $\beta$ are constant, then this construction reduces to the classical algebra of generalized quaternions.

\section{Monogenic Functions in the Commutative Algebra $\mathcal{A}_2(D,\gamma)$}

Next, we generalize elements of the theory of monogenic functions developed in the monograph \cite{Plaksa-Shpakivskyi-2023}
 to the case of a commutative associative algebra with variable structural constants.

We demonstrate our approach for a two-dimensional commutative associative algebra with variable structural constants over the field of real numbers $\mathbb{R}$ and with multiplication table
(\ref{tabl-mnozh-2-vym-algebra-zi-zminnymy-funkciamy}). We denote this algebra (family of algebras) by $\mathcal{A}_2(D,\gamma)$. In what follows, for simplicity of presentation and exposition, we assume that $D\subset\mathbb{R}^2$ and consider the variable
$\zeta=xe_1+ye_2$, $(x,y)\in D$, where $\{e_1,e_2\}$ is the common fixed basis of all algebras in the family $\mathcal{A}_2(D,\gamma)$, with
$\gamma=\{a(x,y),b(x,y)\}$.
Applications may, however, require consideration of other variables, for example,
$\zeta=xi_1+yi_2$, where $i_1,i_2$ are certain vectors of the algebra $\mathcal{A}_2(D,\gamma)$ satisfying prescribed relations.

Consider a function $\Phi:D\to \mathcal{A}_2(D,\gamma)$ which assigns to each point $(x,y)\in D$ an element of the algebra $A_{(x,y)}(\gamma)$, that is,
$\Phi(x,y)\in A_{(x,y)}(\gamma)$.
Since all algebras in the family $\mathcal{A}_2(D,\gamma)$ have the common basis $\{e_1,e_2\}$, the function $\Phi$ can be represented in the form
\begin{equation}\label{rozklad za basysom -2-vym-algebra-zi-zminnymy-funkciamy}
\Phi(\zeta)=U(x,y)e_1+V(x,y)e_2,
\end{equation}
where $U,V:D\to\mathbb{R}$.
The functions $U$ and $V$ will be called the components of the function $\Phi$.

\vskip2mm
\textbf{Definition 1.} Let the components $U,V$ of the function (\ref{rozklad za basysom -2-vym-algebra-zi-zminnymy-funkciamy})
have continuous first-order partial derivatives in the domain $D$. The function $\Phi$ will be called monogenic in $D$
if, at every point $(x,y)\in D$, there exists an element $\Phi'(\zeta)\in A_{(x,y)}(\gamma)$ such that
\begin{equation}\label{ozn-monog-A-2-zmin-koef}
d\Phi=\Phi'(\zeta)d\zeta.
\end{equation}
The element $\Phi'(\zeta)$ will be called the derivative of the function $\Phi$ at the point $\zeta$.

We have
$$
d\Phi=\Phi_xdx+\Phi_ydy, \qquad
d\zeta=dx e_1+dy e_2.
$$
Therefore, from equality (\ref{ozn-monog-A-2-zmin-koef}), we obtain
$$
\Phi_xdx+\Phi_ydy =\Phi'(\zeta)dx e_1+\Phi'(\zeta)dy e_2.
$$
Since $e_1$ is the identity element of the algebra, comparison of the coefficients of $dx$ gives
$\Phi'(\zeta)=\Phi_x$.
Substituting this expression into the equality of the coefficients of $dy$, we obtain
$\Phi_y=\Phi_xe_2$.

Thus, we have the following theorem.

\vskip2mm
\textbf{Theorem 1.} A function $\Phi$ of the form (\ref{rozklad za basysom -2-vym-algebra-zi-zminnymy-funkciamy}) is monogenic in the domain $D$ if and only if the relation
\begin{equation}\label{A-2-zmin-koef-umovy-K-R}
\Phi_y=\Phi_xe_2
\end{equation}
holds.
Conditions (\ref{A-2-zmin-koef-umovy-K-R}) are analogues of the Cauchy--Riemann conditions.

Let us write the Cauchy--Riemann conditions in terms of the components $U,V$ of a monogenic function.

Thus,
$$
U_ye_1+V_ye_2=a(x,y)V_xe_1+\bigl(U_x+b(x,y)V_x\bigr)e_2.
$$
Hence, equating the coefficients of $e_1$ and $e_2$, we obtain the system
\begin{equation}\label{A-2-zmin-koef-systema-K-R}
U_y=a(x,y)V_x,\qquad V_y=U_x+b(x,y)V_x.
\end{equation}

It is clear that the Cauchy--Riemann conditions (\ref{A-2-zmin-koef-umovy-K-R}) or (\ref{A-2-zmin-koef-systema-K-R}) depend on the point $(x,y)\in D$. Thus, in general, the Cauchy--Riemann conditions are different at different points of $D$.

\vskip2mm
\textbf{Example 5.} Consider the function $\Phi(\zeta)=e_2^2$.
According to the multiplication table of the algebra $\mathcal{A}_2(D,\gamma)$, we have
$$
e_2^2=a(x,y)e_1+b(x,y)e_2.
$$

Hence, $\Phi(\zeta)=a(x,y)e_1+b(x,y)e_2$, that is,
$$
U(x,y)=a(x,y),
\qquad
V(x,y)=b(x,y).
$$
Substituting these functions into the generalized Cauchy--Riemann conditions
(\ref{A-2-zmin-koef-systema-K-R}),
we obtain the system
\begin{equation}\label{A-2-zmin-koef-umova-symisnosti}
a_y=ab_x,\qquad
b_y=a_x+bb_x.
\end{equation}

Thus, the function $\Phi(\zeta)=e_2^2$
is monogenic in the domain $D$ if and only if the structural functions $a(x,y)$ and $b(x,y)$ satisfy system (\ref{A-2-zmin-koef-umova-symisnosti}).

Thus, even such a simple function as $e_2^2$ is not always monogenic; its monogenicity requires certain relations between the structural functions.

\vskip2mm
\textbf{Definition 2.} System (\ref{A-2-zmin-koef-umova-symisnosti}) will play a key role throughout the subsequent theory. We call
(\ref{A-2-zmin-koef-umova-symisnosti}) the \textbf{compatibility system} for the structural functions.
The algebra $A_{(x,y)}(\gamma)$ will be called \textbf{compatible} if its structural functions $a(x,y)$ and $b(x,y)$ satisfy the compatibility system (\ref{A-2-zmin-koef-umova-symisnosti}).
The family of all compatible algebras will be denoted by $\mathcal{A}_2^0(D,\gamma)$.

In what follows, we consider only the family $\mathcal{A}_2^0(D,\gamma)$.

We give some additional information about the compatibility system (\ref{A-2-zmin-koef-umova-symisnosti}).
In the case where $b^2+4a\geq0$ in $D$, set
$$
p=\frac{b+\sqrt{b^2+4a}}{2},\qquad
q=\frac{b-\sqrt{b^2+4a}}{2}.
$$
Then
$$
a=-p\,q, \qquad b=p+q.
$$

In these variables, the compatibility system (\ref{A-2-zmin-koef-umova-symisnosti}) takes the form
$$
p_y=pp_x, \qquad
q_y=qq_x.
$$

Thus, the compatibility system splits into two independent Hopf equations (or inviscid Burgers equations). In particular, its general solution can be written implicitly as
$$
p=F(x+py), \qquad q=G(x+qy),
$$
where $F$ and $G$ are arbitrary sufficiently smooth functions.

Let us consider some properties of monogenic functions.

\vskip2mm
\textbf{Theorem 2.}\label{teorema pro sumu monohen funkcij-zmin-koef} Let the functions $\Phi$ and $\Psi$ be monogenic in the domain $D$ in a compatible algebra $\mathcal{A}_2^0(D,\gamma)$.
Then, for arbitrary constants $\alpha,\beta\in\mathbb{R}$, the function
$$
\Theta=\alpha\Phi+\beta\Psi
$$
is also monogenic in the domain $D$.

\textbf{Proof.} Since the functions $\Phi$ and $\Psi$ are monogenic in $D$, they satisfy the Cauchy--Riemann conditions
$$
\Phi_y=\Phi_xe_2,
\qquad
\Psi_y=\Psi_xe_2.
$$

Consider the function
$$
\Theta=\alpha\Phi+\beta\Psi,
\qquad \alpha,\beta\in\mathbb{R}.
$$
Then
$$
\Theta_y=(\alpha\Phi+\beta\Psi)_y
=\alpha\Phi_y+\beta\Psi_y.
$$
Using the monogenicity of $\Phi$ and $\Psi$, we obtain
$$
\Theta_y=\alpha\Phi_xe_2+\beta\Psi_xe_2.
$$

Clearly,
$$
\alpha\Phi_xe_2+\beta\Psi_xe_2=(\alpha\Phi_x+\beta\Psi_x)e_2.
$$
But
$$
\Theta_x=(\alpha\Phi+\beta\Psi)_x=\alpha\Phi_x+\beta\Psi_x.
$$

Hence, $\Theta_y=\Theta_xe_2$.
The theorem is proved.

Let us consider some further examples of monogenic functions.

\vskip2mm
\textbf{Example 6.}
The function $\Phi(\zeta)=\zeta=xe_1+ye_2$ is monogenic in $\mathbb{R}^2$.
Indeed, for the function $\Phi(\zeta)=\zeta$, we have $U(x,y)=x$, $V(x,y)=y$.
Then $U_x=1$, $U_y=0$, $V_x=0$, $V_y=1$.
The Cauchy--Riemann system (\ref{A-2-zmin-koef-systema-K-R})
reduces to the identities $0=a\cdot0$, $1=1+b\cdot0$. Moreover,
$(\zeta)'=e_1=1$.

\vskip2mm
\textbf{Example 7.}\label{pryklad-7-zmin-koef}
Consider the function
$\Phi(\zeta)=\zeta^2=(x^2+ay^2)e_1+(2xy+by^2)e_2$.
Substituting $U(x,y)=x^2+ay^2$ and $V(x,y)=2xy+by^2$ into the Cauchy--Riemann system (\ref{A-2-zmin-koef-systema-K-R}),
we obtain the compatibility system (\ref{A-2-zmin-koef-umova-symisnosti}).
Hence, in every compatible algebra $\mathcal{A}_2^0(D,\gamma)$, the function $\Phi(\zeta)=\zeta^2$ is monogenic in $D$.

\vskip2mm
\textbf{Theorem 3.}\label{teorema pro dobutok monohen funkcij-zmin-koef}
If the functions $\Phi$ and $\Psi$ are monogenic in the domain $D$ in a compatible algebra $\mathcal{A}_2^0(D,\gamma)$, then their product
$\Phi\Psi$ is also monogenic in $D$. Moreover,
$$(\Phi\Psi)'=\Phi'\Psi+\Phi\Psi'.$$

\textbf{Proof.} Let
$$
\Phi(\zeta)=U(x,y)e_1+V(x,y)e_2,
\qquad
\Psi(\zeta)=P(x,y)e_1+Q(x,y)e_2
$$
and
$$
\Theta(\zeta):=\Phi(\zeta)\Psi(\zeta).
$$

Since $e_2^2=a(x,y)e_1+b(x,y)e_2$, we have
$$
\Theta=\Phi\Psi=We_1+Ze_2,
$$
where
$$
W=UP+aVQ,\qquad Z=UQ+VP+bVQ.
$$
Since the functions $\Phi$ and $\Psi$ are monogenic, their components satisfy the Cauchy--Riemann systems
$$
U_y=aV_x,\qquad V_y=U_x+bV_x,
$$
$$
P_y=aQ_x,\qquad Q_y=P_x+bQ_x.
$$

We show that the function $\Theta$ also satisfies the Cauchy--Riemann conditions, namely
$$
W_y=aZ_x,\qquad
Z_y=W_x+bZ_x.
$$
We have
$$
W_y=U_yP+UP_y+a_yVQ+aV_yQ+aVQ_y.
$$
Using the Cauchy--Riemann conditions for $\Phi$ and $\Psi$, we obtain
$$
W_y= aV_xP+aUQ_x+a_yVQ+a(U_x+bV_x)Q+aV(P_x+bQ_x).
$$

On the other hand,
$$
Z_x=U_xQ+UQ_x+V_xP+VP_x+b_xVQ+bV_xQ+bVQ_x.
$$
Hence,
$$
W_y-aZ_x=(a_y-ab_x)VQ.
$$

Since in a compatible algebra $\mathcal{A}_2^0(D,\gamma)$ the structural functions satisfy the compatibility condition
$a_y=ab_x$, we obtain $W_y=aZ_x$.

Now let us verify the second condition. We have
$$
Z_y=U_yQ+UQ_y+V_yP+VP_y+b_yVQ+bV_yQ+bVQ_y.
$$
Substituting the Cauchy--Riemann conditions, we obtain
$$
Z_y=aV_xQ+U(P_x+bQ_x)+(U_x+bV_x)P+aVQ_x+b_yVQ+b(U_x+bV_x)Q+bV(P_x+bQ_x).
$$

On the other hand,
$$
W_x=U_xP+UP_x+a_xVQ+aV_xQ+aVQ_x,
$$
and
$$
Z_x=U_xQ+UQ_x+V_xP+VP_x+b_xVQ+bV_xQ+bVQ_x.
$$
Therefore,
$$
W_x+bZ_x=U_xP+UP_x+a_xVQ+aV_xQ+aVQ_x+bU_xQ+bUQ_x
$$
$$
+bV_xP+bVP_x+bb_xVQ+b^2V_xQ+b^2VQ_x.
$$

Comparing this expression with the formula for $Z_y$, we obtain
$$
Z_y-(W_x+bZ_x)=(b_y-a_x-bb_x)VQ.
$$
Since the structural functions of a compatible algebra also satisfy
$b_y=a_x+bb_x$, it follows that $Z_y=W_x+bZ_x$.

Thus, the components $W$ and $Z$ of the function $\Theta=\Phi\Psi$ satisfy the Cauchy--Riemann system. Therefore, $\Theta$ is monogenic in the domain $D$.

Let us prove the Leibniz formula for the derivative of the product. We have
$\Phi'=\Phi_x$, $\Psi'=\Psi_x$. Therefore,
$$
(\Phi\Psi)'=(\Phi\Psi)_x=\Phi_x\Psi+\Phi\Psi_x=\Phi'\Psi+\Phi\Psi'.
$$
The theorem is proved.

The following corollary follows from the preceding two theorems.

\vskip2mm
\textbf{Corollary 1.}\label{naslidok-1-A-2-zmin}
The set of all functions monogenic in the domain $D$ with values in the compatible algebra $\mathcal{A}_2^0(D,\gamma)$
forms a commutative associative algebra over the field $\mathbb{R}$ with respect to addition of functions, multiplication by real scalars, and pointwise multiplication.

Corollary 1
%\ref{naslidok-1-A-2-zmin}
provides a method for constructing monogenic functions.

It follows from Example 7
%\ref{pryklad-7-zmin-koef}
and Theorem 3
%\ref{teorema pro dobutok monohen funkcij-zmin-koef}
that, for every natural number $n$, the power function
$\zeta^n$ is monogenic in $D$.

Now Theorems 2 and 3
% \ref{teorema pro sumu monohen funkcij-zmin-koef}, \ref{teorema pro dobutok monohen funkcij-zmin-koef}
imply that the polynomial
$$
P_n(\zeta)=\sum\limits_{k=0}^nc_k\zeta^k, \qquad c_k\in\mathbb{R}
$$
in the compatible algebra $\mathcal{A}_2^0(D,\gamma)$ is monogenic in $D$. Moreover,
$$
(P_n(\zeta))'=\sum\limits_{k=1}^nkc_k\zeta^{k-1}.
$$

We show that the compatibility condition (\ref{A-2-zmin-koef-umova-symisnosti}) is not only sufficient but also necessary for the product of two monogenic functions to be monogenic again.

\vskip2mm
\textbf{Theorem 4.}\label{teorema pro dobutok monohen funkcij-zmin-koef-obernena}
Let two monogenic functions $\Phi,\Psi$ be given in the two-dimensional commutative algebra with variable structural constants $\mathcal{A}_2(D,\gamma)$.
Then the product $\Phi\Psi$ of any two monogenic functions is monogenic if and only if the compatibility condition (\ref{A-2-zmin-koef-umova-symisnosti}) holds, that is,
$\mathcal{A}_2(D,\gamma)\equiv\mathcal{A}^0_2(D,\gamma)$.

\textbf{Proof.} \textit{Sufficiency} was proved in Theorem 3.
%\ref{teorema pro dobutok monohen funkcij-zmin-koef}.

We prove \textit{necessity}.
Suppose that the product of any two monogenic functions is again monogenic. Take arbitrary monogenic functions
$$
\Phi=Ue_1+Ve_2,\qquad \Psi=Se_1+Te_2.
$$
Then their components satisfy the Cauchy--Riemann conditions
$$
U_y=aV_x,\qquad V_y=U_x+bV_x,
$$
$$
S_y=aT_x,\qquad T_y=S_x+bT_x.
$$

Their product has the form
$$
\Phi\Psi=(US+aVT)e_1+(UT+VS+bVT)e_2.
$$
Set
$$
P=US+aVT,\qquad Q=UT+VS+bVT.
$$
By assumption, the function $\Phi\Psi=Pe_1+Qe_2$ is monogenic, and therefore
$$
P_y=aQ_x,\qquad Q_y=P_x+bQ_x.
$$

On the other hand, a direct computation using the monogenicity conditions for $\Phi$ and $\Psi$ gives
$$
P_y-aQ_x=(a_y-a\,b_x)VT,
$$
as well as
$$
Q_y-P_x-bQ_x=(b_y-a_x-b\,b_x)VT.
$$
Since the function $\Phi\Psi$ is monogenic, the left-hand sides of these equalities are equal to zero. Hence,
$$
(a_y-a\,b_x)VT=0,\qquad
(b_y-a_x-b\,b_x)VT=0.
$$

Since monogenic functions can be chosen so that $VT\not\equiv0$, the compatibility conditions (\ref{A-2-zmin-koef-umova-symisnosti}) necessarily follow.
Thus, the compatibility condition is necessary.

This is precisely why the compatibility conditions (\ref{A-2-zmin-koef-umova-symisnosti}) play a key role in our investigations.

\vskip2mm
\textbf{Theorem 5.}
Let the structural constants $a$ and $b$ of the compatible algebra $\mathcal{A}_2^0(D,\gamma)$ be constant.
If the function $\Phi(x,y)=u(x,y)e_1+v(x,y)e_2$
is monogenic and its components $u,v$ possess continuous partial derivatives up to and including second order, then its derivative
$$
\Phi'(Z)=\Phi_x(x,y)=u_x(x,y)e_1+v_x(x,y)e_2
$$
is also monogenic.

\textbf{Proof.} Since $\Phi$ is monogenic, its components satisfy the Cauchy--Riemann conditions (\ref{A-2-zmin-koef-systema-K-R}):
$$
u_y=av_x,\qquad v_y=u_x+bv_x.
$$
Differentiating these equalities with respect to $x$, we obtain
$$
u_{xy}=av_{xx}, \qquad
v_{xy}=u_{xx}+bv_{xx}.
$$
But these are precisely the monogenicity conditions for the function
$\Phi_x=u_xe_1+v_xe_2$.
The theorem is proved.

\vskip2mm
\textbf{Remark 2.}
If the structural constants $a$ and $b$ are functions, then the preceding theorem is, in general, no longer valid.

\section{Relationship Between Monogenic Functions and Linear PDEs with Variable Coefficients}

In \cite{Pogorui-2017}, the hypercomplex method was extended to certain partial differential equations with variable coefficients, 
and their solutions were constructed using monogenic functions in commutative algebras. In paper \cite{Joseph}, 
several families of exact solutions are constructed for second-order partial differential equations with variable coefficients.
In paper \cite{Reutskiy}, the method of particular solutions is proposed for solving second- and fourth-order partial differential equations with variable coefficients.

\vskip2mm

In this section, we establish a relationship between monogenic functions in the algebra
$\mathcal{A}^0_2(D,\gamma)$ and linear PDEs with variable coefficients. It is shown that the components of a monogenic
function satisfy certain linear PDEs with variable coefficients. This makes it possible to construct explicit
solutions of these classes of PDEs by means of the theory of monogenic functions.

\vskip2mm
\textbf{Theorem 6.}\label{zvyazok z rivnianniam-A-2-0-zmin-funk} 
If the function $\Phi(x,y)=U(x,y)e_1+V(x,y)e_2$
is monogenic in the compatible algebra $\mathcal{A}^0_2(D,\gamma)$, then its component $V$ satisfies the equation
\begin{equation}\label{V-zmin-funk}
V_{yy}=aV_{xx}+bV_{xy}+\left(2a_x+b\,b_x\right)V_x,
\end{equation}
whereas the component $U$ satisfies the equation
\begin{equation}\label{U-zmin-funk}
aU_{yy}=a^2U_{xx}+abU_{xy}+\left(2ab_x-ba_x\right)U_y.
\end{equation}

\textbf{Proof.} Since the function $\Phi=Ue_1+Ve_2$ is monogenic, its components satisfy the Cauchy--Riemann conditions
$$
U_y=aV_x,\qquad V_y=U_x+bV_x.
$$

First, we derive the equation for $V$. Differentiating the second equality with respect to $y$, we obtain
$$
V_{yy}=U_{xy}+b_yV_x+bV_{xy}.
$$
From the first equality, we have
$$
U_{xy}=(aV_x)_x=a_xV_x+aV_{xx}.
$$
Therefore,
$$
V_{yy}=aV_{xx}+bV_{xy}+(a_x+b_y)V_x.
$$
Using the compatibility condition $b_y=a_x+b\,b_x$,
we obtain
$$
V_{yy}=aV_{xx}+bV_{xy}+\left(2a_x+b\,b_x\right)V_x.
$$

Now we derive the equation for $U$. Differentiating the first equality with respect to $y$, we obtain
$$
U_{yy}=a_yV_x+aV_{xy}.
$$
Differentiating the second equality with respect to $x$, we get
$$
V_{xy}=U_{xx}+b_xV_x+bV_{xx}.
$$
Moreover, from the equality $U_y=aV_x$, we have
$$
U_{xy}=a_xV_x+aV_{xx}.
$$
Hence,
$$
aV_{xx}=U_{xy}-a_xV_x.
$$
Therefore,
$$
U_{yy}=a_yV_x+a\left(U_{xx}+b_xV_x+bV_{xx}\right)
=aU_{xx}+bU_{xy}+\left(a_y+ab_x-ba_x\right)V_x.
$$
Since, by the compatibility condition, $a_y=a\,b_x$,
we obtain
$$
U_{yy}=aU_{xx}+bU_{xy}+\left(2ab_x-ba_x\right)V_x.
$$
Multiplying this equality by $a$ and using $U_y=aV_x$, we obtain
$$
aU_{yy}=a^2U_{xx}+abU_{xy}+\left(2ab_x-ba_x\right)U_y.
$$
The theorem is proved.

We now present converse results.

\vskip2mm
\textbf{Theorem 7.} Let functions $a(x,y)$ and $b(x,y)$ be given in a simply connected domain $D\subset\mathbb{R}^2$, with $a(x,y)\neq 0$ in $D$, and suppose that the compatibility conditions
(\ref{A-2-zmin-koef-umova-symisnosti}) hold. Let, in addition, the function $U=U(x,y)\in C^2(D)$ satisfy equation (\ref{U-zmin-funk}).
Then there exists a function $V=V(x,y)$ such that $\Phi(x,y)=U(x,y)e_1+V(x,y)e_2$ is monogenic in $D$.

\textbf{Proof.} In order for the function $\Phi=Ue_1+Ve_2$
to be monogenic, its components must satisfy the Cauchy--Riemann system
$$
U_y=aV_x,\qquad V_y=U_x+bV_x.
$$
Since $a\neq0$, the first equality gives
$V_x=\frac{U_y}{a}$.
Then the second equality takes the form
$V_y=U_x+\frac{b}{a}U_y$.
Thus, the function $V$ must be found from the system
$$
V_x=\frac{U_y}{a},\qquad
V_y=U_x+\frac{b}{a}U_y.
$$

Consider the differential form
$$
\omega=\frac{U_y}{a}\,dx+\left(U_x+\frac{b}{a}U_y\right)dy.
$$
The function $V$ exists if and only if the form $\omega$ is exact. Since the domain $D$ is simply connected, it is sufficient to verify that $\omega$ is closed:
$$
\left(\frac{U_y}{a}\right)_y= \left(U_x+\frac{b}{a}U_y\right)_x.
$$

We compute
$$
\left(\frac{U_y}{a}\right)_y=\frac{U_{yy}}{a}-\frac{a_y}{a^2}U_y,
$$
and also
$$
\left(U_x+\frac{b}{a}U_y\right)_x
=U_{xx}+\frac{b}{a}U_{xy}
+\left(\frac{b}{a}\right)_xU_y.
$$
Therefore, the closedness condition has the form
$$
\frac{U_{yy}}{a}-\frac{a_y}{a^2}U_y
=U_{xx}+\frac{b}{a}U_{xy}
+\left(\frac{b}{a}\right)_xU_y.
$$
Multiplying this equality by $a^2$, we obtain
$$
aU_{yy}-a_yU_y
=a^2U_{xx}+abU_{xy}+(ab_x-ba_x)U_y.
$$
Hence,
$$
aU_{yy}= a^2U_{xx}+abU_{xy}+(a_y+ab_x-ba_x)U_y.
$$
Since, by the compatibility condition, $a_y=ab_x$, this equality becomes
$$
aU_{yy}= a^2U_{xx}+abU_{xy}+\left(2ab_x-ba_x\right)U_y.
$$
This is precisely the equation which, by assumption, is satisfied by the function $U$.

Hence, the form $\omega$ is closed, and therefore, in the simply connected domain $D$, there exists a function $V$ such that $dV=\omega$.
That is,
$$
V_x=\frac{U_y}{a},\qquad
V_y=U_x+\frac{b}{a}U_y.
$$
Hence,
$$
U_y=aV_x,\qquad V_y=U_x+bV_x.
$$
Therefore, the function $\Phi=Ue_1+Ve_2$ is monogenic.
The theorem is proved.

Moreover, one can write an explicit formula for computing the function $V$:
$$
V(x,y)=V(x_0,y_0)+ \int_{\gamma}\frac{U_y}{a}\,dx+
\left(U_x+\frac{b}{a}U_y\right)dy,
$$
where $\gamma$ is an arbitrary curve in $D$ joining $(x_0,y_0)$ to $(x,y)$. Since the form is closed, the integral is independent of the choice of path.

\vskip2mm
\textbf{Theorem 8.} Let functions $a(x,y)$ and $b(x,y)$ be given in a simply connected domain $D$ and satisfy the compatibility conditions
(\ref{A-2-zmin-koef-umova-symisnosti}).
Let, in addition, the function $V=V(x,y)\in C^2(D)$ satisfy equation (\ref{V-zmin-funk}).
Then there exists a function $U=U(x,y)$ such that $\Phi(x,y)=U(x,y)e_1+V(x,y)e_2$
is monogenic in $D$.

\textbf{Proof.} In order for the function $\Phi=Ue_1+Ve_2$
to be monogenic, its components must satisfy the Cauchy--Riemann system
$$
U_y=aV_x,\qquad V_y=U_x+bV_x.
$$
From the second equality, we obtain $U_x=V_y-bV_x$.
Thus, the function $U$ must be found from the system
$$
U_x=V_y-bV_x,\qquad U_y=aV_x.
$$

Consider the differential form
$$
\omega=\left(V_y-bV_x\right)dx+aV_x\,dy.
$$
The function $U$ exists if and only if the form $\omega$ is exact. Since the domain $D$ is simply connected, it is sufficient to verify its closedness:
$$
\left(V_y-bV_x\right)_y=(aV_x)_x.
$$

We compute
$$
\left(V_y-bV_x\right)_y=V_{yy}-b_yV_x-bV_{xy},
$$
and also
$$
(aV_x)_x=a_xV_x+aV_{xx}.
$$
Therefore, the closedness condition has the form
$$
V_{yy}-b_yV_x-bV_{xy}=a_xV_x+aV_{xx}.
$$
Hence,
$$
V_{yy}=aV_{xx}+bV_{xy}+(a_x+b_y)V_x.
$$
Since, by the compatibility condition, $b_y=a_x+b\,b_x$,
we obtain
$$
V_{yy}=aV_{xx}+bV_{xy}+\left(2a_x+b\,b_x\right)V_x.
$$
This is precisely the equation which, by assumption, is satisfied by the function $V$.

Hence, the form $\omega$ is closed, and therefore, in the simply connected domain $D$, there exists a function $U$ such that
$dU=\omega$.
That is,
$$
U_x=V_y-bV_x,\qquad U_y=aV_x.
$$
Hence,
$$
U_y=aV_x,\qquad V_y=U_x+bV_x.
$$
Therefore, the function $\Phi=Ue_1+Ve_2$ is monogenic.
The theorem is proved.

The function $U$ can be computed by the formula
$$
U(x,y)=U(x_0,y_0)+\int_{\gamma} \left(V_y-bV_x\right)dx+aV_x\,dy,
$$
where $\gamma$ is an arbitrary curve in $D$ joining $(x_0,y_0)$ to $(x,y)$. Since the form is closed, the integral is independent of the choice of path.

\vskip2mm
\textbf{Remark 5.} Thus, under the condition $a(x,y)\neq0$ in $D$, the set of twice continuously differentiable solutions of equations (\ref{U-zmin-funk}) and (\ref{V-zmin-funk})
in a simply connected domain $D$ is in one-to-one correspondence with the set of
monogenic functions in the compatible algebra $\mathcal{A}^0_2(D,\gamma)$.
This means that all solutions of equations (\ref{U-zmin-funk}) and (\ref{V-zmin-funk}) in the class $C^2(D)$ can be obtained as components of monogenic functions
$\Phi$ with values in the compatible algebra $\mathcal{A}^0_2(D,\gamma)$.
Corollary 1
%\ref{naslidok-1-A-2-zmin}
provides a method for constructing monogenic functions.

We next give a theorem on properties of monogenic functions in the compatible algebra $\mathcal{A}^0_2(D,\gamma)$, which is an analogue of the Cauchy integral theorem.

\vskip2mm
\textbf{Definition 3.} Let $\gamma$ be a piecewise smooth curve in $\mathbb{R}^2$. The integral of the function $\Phi(\zeta)=U(x,y)e_1+V(x,y)e_2$ along the curve $\gamma$
is defined by
\begin{equation}\label{oznachen integrala po kryviy-A-2-0}
\int\limits_{\gamma}\Phi\,d\zeta
=e_1\int\limits_{\gamma}U\,dx+aV\,dy
+e_2\int\limits_{\gamma}V\,dx+(U+bV)\,dy,
\end{equation}
provided that the integrals on the right-hand side exist.

\vskip2mm
\textbf{Theorem 9.} Let the function $\Phi(x,y)=U(x,y)e_1+V(x,y)e_2$
be monogenic in a domain $D$ in the compatible algebra $\mathcal{A}_2^0(D,\gamma)$, and let $\gamma$ be a positively oriented closed piecewise smooth curve in $D$
bounding a domain $G\subset D$.
Then
\begin{equation}\label{int-teor-koshi-zmin-func}
\int\limits_{\gamma}\Phi\,d\zeta
=\int\hspace{-3mm}\int\limits_{G} V(a_x e_1+b_xe_2)\,dxdy.
\end{equation}

\textbf{Proof.} By the definition of the integral, we have
$$
\int\limits_{\gamma}\Phi\,d\zeta=
e_1\int\limits_{\gamma}U\,dx+aV\,dy
+ e_2\int\limits_{\gamma}V\,dx+(U+bV)\,dy.
$$
Applying Green's formula to the first integral, we obtain
$$
\int\limits_{\gamma}U\,dx+aV\,dy=
\int\hspace{-3mm}\int\limits_G\biggr((aV)_x-U_y\biggr)\,dxdy.
$$
Since $(aV)_x=a_xV+aV_x$ and, by monogenicity, $U_y=aV_x$,
we have
$$
(aV)_x-U_y=a_xV.
$$
Therefore,
$$
\int\limits_{\gamma}U\,dx+aV\,dy
=
\int\hspace{-3mm}\int\limits_G a_xV\,dxdy.
$$

Similarly,
$$
\int\limits_{\gamma}V\,dx+(U+bV)\,dy
=
\int\hspace{-3mm}\int\limits_G\biggr((U+bV)_x-V_y\biggr)\,dxdy.
$$
But
$$
(U+bV)_x=U_x+b_xV+bV_x,
$$
whereas monogenicity gives
$V_y=U_x+bV_x$.
Hence,
$(U+bV)_x-V_y=b_xV$.
Therefore,
$$
\int\limits_{\gamma}V\,dx+(U+bV)\,dy
=
\int\hspace{-3mm}\int\limits_G b_xV\,dxdy.
$$
Substituting these equalities into the definition of the integral, we obtain (\ref{int-teor-koshi-zmin-func}).
The theorem is proved.

\section{Examples}
%\label{Paragraf-6}

We give several examples of structural functions $a$ and $b$ satisfying the compatibility system (\ref{A-2-zmin-koef-umova-symisnosti}) and the corresponding equations for the components of monogenic functions.

\vskip2mm
\textbf{Example 8.}

Let $a=1$, $b=0$. Then the equations for the components of a monogenic function take the form
$$
U_{yy}=U_{xx}, \qquad V_{yy}=V_{xx}.
$$
Thus, both components of a monogenic function are solutions of the one-dimensional wave equation. In this case, the family of algebras
(\ref{tabl-mnozh-2-vym-algebra-zi-zminnymy-funkciamy}) reduces to the algebra of double numbers. Hence, all solutions of the wave equation
can be obtained as components of monogenic functions in the algebra of double numbers.

\vskip2mm
\textbf{Example 9.} Let
$$
p=\frac{x}{1-y},\qquad q=1.
$$
Then
$$
p_y=pp_x,\qquad q_y=qq_x=0,
$$
and therefore the functions
$$
a=-pq=-\frac{x}{1-y},\qquad b=p+q=\frac{x}{1-y}+1
$$
satisfy the compatibility system (\ref{A-2-zmin-koef-umova-symisnosti}). Hence, in this case, the family of algebras
(\ref{tabl-mnozh-2-vym-algebra-zi-zminnymy-funkciamy}) takes the form
\begin{equation}\label{tabl-mnozh-pryklad-9-zmin-funk}
\begin{array}{c||c|c|}
\cdot & e_1 & e_2\\
\hline\hline
e_1 & e_1 & e_2\\
\hline
e_2 & e_2 & -\frac{x}{1-y}\,e_1+\left(\frac{x}{1-y}+1\right)e_2\\
\hline
\end{array}\,.
\end{equation}

For the derivatives, we have
$$
a_x=-\frac1{1-y},\qquad b_x=\frac1{1-y}.
$$
Therefore, the equation for the component $V$ takes the form
\begin{equation}\label{tabl-mnozh-pryklad-9-zmin-funk-V}
V_{yy}= -\frac{x}{1-y}V_{xx}+ \left(1+\frac{x}{1-y}\right)V_{xy}+ \frac{x+y-1}{(1-y)^2}V_x.
\end{equation}

The component $U$ satisfies the equation
\begin{equation}\label{tabl-mnozh-pryklad-9-zmin-funk-U}
-\frac{x}{1-y}U_{yy}= \frac{x^2}{(1-y)^2}U_{xx}- \frac{x}{1-y}\left(1+\frac{x}{1-y}\right)U_{xy}+\frac{1-y-x}{(1-y)^2}U_y.
\end{equation}

Let us give examples of some solutions of equations (\ref{tabl-mnozh-pryklad-9-zmin-funk-V}) and (\ref{tabl-mnozh-pryklad-9-zmin-funk-U}).
To this end, in the algebra (\ref{tabl-mnozh-pryklad-9-zmin-funk}) consider the monogenic function
\begin{equation}\label{tabl-mnozh-pryklad-9-zeta-2}
\Phi(\zeta)=\zeta^2=(x^2+ay^2)e_1+(2xy+by^2)e_2,
\end{equation}
that is,
\begin{equation}\label{tabl-mnozh-pryklad-9-zmin-funk-U-V}
U=x^2-\frac{xy^2}{1-y},\qquad V=2xy+\left(\frac{x}{1-y}+1\right)y^2.
\end{equation}

Let us compute the derivatives of the component $V$:
$$
V_x=\frac{y(2-y)}{1-y},\quad V_{xx}=0, \quad V_{xy}=\frac{y^2-2y+2}{(1-y)^2},\quad V_{yy}=2+\frac{2x}{(1-y)^3}.
$$
Then the right-hand side of equation (\ref{tabl-mnozh-pryklad-9-zmin-funk-V}) is equal to
$$
-\frac{x}{1-y}\cdot 0+\left(1+\frac{x}{1-y}\right)\frac{y^2-2y+2}{(1-y)^2}+\frac{x+y-1}{(1-y)^2}\frac{y(2-y)}{1-y}.
$$
After reducing to a common denominator, we obtain
$$
2+\frac{2x}{(1-y)^3}.
$$
Hence, the function $V$ from (\ref{tabl-mnozh-pryklad-9-zmin-funk-U-V}) is indeed a solution of equation (\ref{tabl-mnozh-pryklad-9-zmin-funk-V}).

Now let us verify the equation for $U$. We have
$$
U_y=\frac{xy(y-2)}{(1-y)^2},\quad U_{xx}=2,\quad U_{xy}=\frac{y(y-2)}{(1-y)^2},\quad U_{yy}=-\frac{2x}{(1-y)^3}.
$$
Therefore, the left-hand side of the equation is
$$
-\frac{x}{1-y}U_{yy}=-\frac{x}{1-y}\left(-\frac{2x}{(1-y)^3}\right)=\frac{2x^2}{(1-y)^4}.
$$

Substituting the calculated partial derivatives into the right-hand side of equation (\ref{tabl-mnozh-pryklad-9-zmin-funk-U}), we obtain
$$
\frac{x^2}{(1-y)^2}\cdot 2-\frac{x}{1-y}\left(1+\frac{x}{1-y}\right)
\frac{y(y-2)}{(1-y)^2}+\frac{1-y-x}{(1-y)^2}\frac{xy(y-2)}{(1-y)^2}.
$$
After simplification, this expression is equal to
$$
\frac{2x^2}{(1-y)^4}.
$$
Hence, the function $U$ from (\ref{tabl-mnozh-pryklad-9-zmin-funk-U-V}) is indeed a solution of equation (\ref{tabl-mnozh-pryklad-9-zmin-funk-U}).

Thus, the components of the function $\Phi(\zeta)=\zeta^2$ indeed satisfy equations
(\ref{tabl-mnozh-pryklad-9-zmin-funk-V}) and (\ref{tabl-mnozh-pryklad-9-zmin-funk-U}).

\vskip2mm
\textbf{Example 10.} Let
$$
p=\frac{x}{1-y},\qquad q=\frac{2x}{1-2y}.
$$
Then
$$
p_y=pp_x,\qquad q_y=qq_x.
$$
Set
$$
a=-pq=-\frac{2x^2}{(1-y)(1-2y)},\qquad
b=p+q=\frac{x}{1-y}+\frac{2x}{1-2y}.
$$
Then the functions $a$ and $b$ satisfy the compatibility system (\ref{A-2-zmin-koef-umova-symisnosti}).

Hence, in this case, the family of algebras
(\ref{tabl-mnozh-2-vym-algebra-zi-zminnymy-funkciamy}) takes the form
\begin{equation}\label{tabl-mnozh-pryklad-10-zmin-funk}
\begin{array}{c||c|c|}
\cdot & e_1 & e_2\\
\hline\hline
e_1 & e_1 & e_2\\
\hline
e_2 & e_2 & -\frac{2x^2}{(1-y)(1-2y)}\,e_1+\left(\frac{x}{1-y}+\frac{2x}{1-2y}\right)e_2\\
\hline
\end{array}\,.
\end{equation}

For the derivatives, we have
$$
a_x=-\frac{4x}{(1-y)(1-2y)},\qquad b_x=\frac{1}{1-y}+\frac{2}{1-2y}.
$$

Therefore, the equation for the component $V$ has the form
\begin{equation}\label{tabl-mnozh-pryklad-10-zmin-funk-V}
V_{yy}=-\frac{2x^2}{(1-y)(1-2y)}V_{xx}+\left(\frac{x}{1-y}+\frac{2x}{1-2y}\right)V_{xy}+\frac{x}{(1-y)^2(1-2y)^2}V_x.
\end{equation}

The equation for the component $U$ has the form
\begin{equation}\label{tabl-mnozh-pryklad-10-zmin-funk-U}
-\frac{2x^2}{(1-y)(1-2y)}U_{yy}=\frac{4x^4}{(1-y)^2(1-2y)^2}U_{xx}+\frac{2x^3(4y-3)}{(1-y)^2(1-2y)^2}U_{xy}.
\end{equation}

Let us give examples of some solutions of equations (\ref{tabl-mnozh-pryklad-10-zmin-funk-V}) and (\ref{tabl-mnozh-pryklad-10-zmin-funk-U}).
To this end, in the algebra (\ref{tabl-mnozh-pryklad-10-zmin-funk}), we again consider the monogenic function
(\ref{tabl-mnozh-pryklad-9-zeta-2}).

In this algebra, the components of the function (\ref{tabl-mnozh-pryklad-9-zeta-2}) have the form
\begin{equation}\label{tabl-mnozh-pryklad-10-zmin-funk-U-V}
U=x^2-\frac{2x^2y^2}{(1-y)(1-2y)},\quad
V=2xy+\left(\frac{x}{1-y}+\frac{2x}{1-2y}\right)y^2.
\end{equation}

Let us compute the derivatives of the component $V$:
$$
V_x=\frac{y(2-3y)}{(1-y)(1-2y)},\quad V_{xx}=0,\quad V_{xy}=\frac{5y^2-6y+2}{(1-y)^2(1-2y)^2},
$$
$$
V_{yy}=-\frac{2x(10y^3-18y^2+12y-3)}{(1-y)^3(1-2y)^3}.
$$

Then the right-hand side of equation (\ref{tabl-mnozh-pryklad-10-zmin-funk-V}) is equal to
$$
-\frac{2x^2}{(1-y)(1-2y)}\cdot 0+\left(\frac{x}{1-y}+\frac{2x}{1-2y}\right)\frac{5y^2-6y+2}{(1-y)^2(1-2y)^2}
$$
$$
+\frac{x}{(1-y)^2(1-2y)^2}\frac{y(2-3y)}{(1-y)(1-2y)}.
$$
After reducing to a common denominator, we obtain
$$
-\frac{2x(10y^3-18y^2+12y-3)}
{(1-y)^3(1-2y)^3}.
$$
Hence, the component $V$ indeed satisfies equation
(\ref{tabl-mnozh-pryklad-10-zmin-funk-V}).

Now let us verify the equation for the component $U$. We have
$$
U_y=\frac{2x^2y(3y-2)}{(1-y)^2(1-2y)^2},\quad
U_{xx}=\frac{2(1-3y)}{(1-y)(1-2y)},
$$
$$
U_{xy}=\frac{4xy(3y-2)}{(1-y)^2(1-2y)^2},\quad U_{yy}=-\frac{4x^2(6y^3-6y^2+1)}{(1-y)^3(1-2y)^3}.
$$

The left-hand side of equation (\ref{tabl-mnozh-pryklad-10-zmin-funk-U}) is
$$
-\frac{2x^2}{(1-y)(1-2y)}U_{yy}
=-\frac{2x^2}{(1-y)(1-2y)}
\left(-\frac{4x^2(6y^3-6y^2+1)}{(1-y)^3(1-2y)^3}\right).
$$
Thus,
$$
-\frac{2x^2}{(1-y)(1-2y)}U_{yy}
=
\frac{8x^4(6y^3-6y^2+1)}{(1-y)^4(1-2y)^4}.
$$

The right-hand side of equation (\ref{tabl-mnozh-pryklad-10-zmin-funk-U}), after substitution of the derivatives, has the form
$$
\frac{4x^4}{(1-y)^2(1-2y)^2}\frac{2(1-3y)}{(1-y)(1-2y)}
+\frac{2x^3(4y-3)}{(1-y)^2(1-2y)^2}\frac{4xy(3y-2)}{(1-y)^2(1-2y)^2}.
$$
After simplification, we obtain
$$
\frac{8x^4(6y^3-6y^2+1)}{(1-y)^4(1-2y)^4}.
$$
Hence, the component $U$ from (\ref{tabl-mnozh-pryklad-10-zmin-funk-U-V}) indeed satisfies equation (\ref{tabl-mnozh-pryklad-10-zmin-funk-U}).

Thus, the components of the function $\Phi(\zeta)=\zeta^2$ satisfy equations
(\ref{tabl-mnozh-pryklad-10-zmin-funk-V}) and (\ref{tabl-mnozh-pryklad-10-zmin-funk-U}).

\vskip2mm
\textbf{Example 11.} Let
$$
p=\frac{\alpha x}{1-\alpha y},\qquad
q=\frac{\beta x}{1-\beta y},
$$
where $\alpha,\beta$ are real constants. Then
$$
p_y=pp_x,\qquad q_y=qq_x.
$$
Therefore, the functions
$$
a=-pq=-\frac{\alpha\beta x^2}{(1-\alpha y)(1-\beta y)},
$$
$$
b=p+q=\frac{\alpha x}{1-\alpha y}+\frac{\beta x}{1-\beta y}
$$
satisfy the compatibility system (\ref{A-2-zmin-koef-umova-symisnosti}).

Knowing the structural functions $a$ and $b$, one can easily write down the corresponding multiplication table for the family of algebras (\ref{tabl-mnozh-2-vym-algebra-zi-zminnymy-funkciamy}).

We have
$$
a_x=-\frac{2\alpha\beta x}{(1-\alpha y)(1-\beta y)},\qquad
b_x=\frac{\alpha}{1-\alpha y}+\frac{\beta}{1-\beta y}.
$$

Hence, the equation for the component $V$ has the form
\begin{equation}\label{tabl-mnozh-pryklad-11-zmin-funk-V}
V_{yy}=-\frac{\alpha\beta x^2}{(1-\alpha y)(1-\beta y)}V_{xx}+
\left(\frac{\alpha x}{1-\alpha y}+\frac{\beta x}{1-\beta y}\right)V_{xy}+
\frac{x(\alpha-\beta)^2}{(1-\alpha y)^2(1-\beta y)^2}V_x.
\end{equation}

The equation for the component $U$ has the form
$$
-\frac{\alpha\beta x^2}{(1-\alpha y)(1-\beta y)}U_{yy}
=
\frac{\alpha^2\beta^2 x^4}{(1-\alpha y)^2(1-\beta y)^2}U_{xx}
$$
\begin{equation}\label{tabl-mnozh-pryklad-11-zmin-funk-U}
-\frac{\alpha\beta x^3}{(1-\alpha y)(1-\beta y)}
\left(\frac{\alpha}{1-\alpha y}+\frac{\beta}{1-\beta y}\right)U_{xy}.
\end{equation}

For this family of algebras, it is not difficult to show that the components $V$ and $U$ of the monogenic function (\ref{tabl-mnozh-pryklad-9-zeta-2})
are solutions of equations (\ref{tabl-mnozh-pryklad-11-zmin-funk-V}) and (\ref{tabl-mnozh-pryklad-11-zmin-funk-U}), respectively.

\section{Concluding Remarks}

1. If we want monogenic functions to be associated, instead of the differential equations (\ref{V-zmin-funk}) and (\ref{U-zmin-funk}), with a differential equation involving a larger number of independent variables $x$, $y$, $z$,\ldots, then it is necessary to consider algebras with variable structural constants of higher dimension.

2. To investigate solutions of the differential equations (\ref{V-zmin-funk}) and (\ref{U-zmin-funk}) or other higher-dimensional differential equations, it becomes necessary to develop an analogue of commutative hypercomplex analysis for the families of algebras introduced above.

\subsection*{Acknowledgment}
 This work was supported by a grant from the Simons Foundation 
(SFI-PD-Ukraine-00014586,V.S.Sh.). The author was supported by budget program “Support for the development of priority areas of research” 
(KPKVK 6541230).

\renewcommand{\refname}{References}

\end{document}